\documentclass[12pt, a4paper, leqno]{article}
\usepackage[intlimits]{amsmath}
\usepackage{amsthm}
\usepackage{amsfonts}
\DeclareMathOperator{\im}{Im}
\DeclareMathOperator{\re}{Re}
\DeclareMathOperator{\Tr}{Tr}
\DeclareMathOperator{\dist}{dist}

\DeclareMathOperator{\diver}{div}
\newcommand{\bgrad}{\text{\bf grad}\,}
\newcommand{\HC}{\mathcal{H}}
\newcommand{\alp}{\alpha}

\newcommand{\lam}{\lambda}

\newcommand{\Del}{\Delta}
\newcommand{\Lam}{\Lambda}
\newcommand{\Ome}{\Omega}
\newcommand{\ubold}{\mathbf{u}}
\newcommand{\xbold}{\mathbf{x}}
\newcommand{\Rbb}[1][]{\mathbb{R}^{#1}}
\newcommand{\dpar}{\partial}
\allowdisplaybreaks[4]
\newcommand{\scal}[1]{\ensuremath{\langle #1\rangle}}
\newcommand{\Om}{\Omega}
\DeclareMathOperator{\spec}{spec}
\numberwithin{equation}{section}
\theoremstyle{plain}
\newtheorem{thm}{Theorem}[section]

\newtheorem{lem}[thm]{Lemma}

\newtheorem{cor}[thm]{Corollary}
\theoremstyle{definition}

\newtheorem{exmp}[thm]{Example}
\theoremstyle{remark}
\newtheorem{rem}[thm]{Remark}

\begin{document}
\title{Commutators, Spectral Trace Identities, and Universal Estimates for Eigenvalues}
\author{Michael Levitin\thanks{The research of M.L. was
partially supported by EPSRC grant GR/M20990}\\
Department of Mathematics, Heriot-Watt University\\
Riccarton, Edinburgh EH14 4AS,
U.~K.
\\
email M.Levitin@ma.hw.ac.uk
\and
Leonid Parnovski\thanks{The research of L.P. was
partially supported by EPSRC grant GR/M20549}\\
Department of Mathematics, University College London\\
Gower Street, London WC1E 6BT, U.~K.\\
email Leonid@math.ucl.ac.uk}
%
%
%
\date{21 February 2001}
\maketitle
\begin{abstract} Using simple commutator relations, we obtain
several trace identities involving eigenvalues and eigenfunctions of an
abstract self-adjoint operator acting in a Hilbert space.
Applications involve abstract universal estimates for
the eigenvalue gaps. As particular examples, we present
simple proofs of the classical universal estimates for
eigenvalues of the Dirichlet Laplacian, as well as of some
known and new results for other differential operators and
systems. We also suggest an extension of the methods to the
case of non-self-adjoint operators.
\end{abstract}
\section{Introduction}

In 1956, Payne, P\'olya and Weinberger \cite{PaPoWe} have shown that if
$\{\lam_j\}$ is the set of (positive) eigenvalues of the Dirichlet boundary
value problem for the Laplacian in a domain $\Om\subset\Rbb[n]$, then
\begin{equation*}\label{eq:PPW}
\lam_{m+1}-\lam_m\le \frac{4}{mn}\sum_{j=1}^m\lam_j\tag{PPW}
\end{equation*}
for each $m=1,2,\dots$.

This inequality was improved to
\begin{equation*}\label{eq:HP}
\sum_{j=1}^m \frac{\lam_j}{\lam_{m+1}-\lam_j}\ge \frac{mn}{4}.\tag{HP}
\end{equation*}
by Hile and Protter \cite{HiPr}. This is indeed stronger than \eqref{eq:PPW},
which is
obtained from \eqref{eq:HP} by replacing all $\lam_j$ in the denominators in the
left-hand side by $\lam_m$.

Later, Hongcang Yang \cite{Ya} proved an even stronger inequality
\begin{equation*}\label{eq:HCY1}
\sum_{j=1}^m \left(\lam_{m+1}-\lam_j\right)\left(\lam_{m+1}-\left(1+\frac{4}{n}\right)
\lam_j\right)\le 0\,,\tag{HCY-1}
\end{equation*}
which after some modifications implies an explicit estimate
\begin{equation*}\label{eq:HCY2}
\lam_{m+1}\le\left(1+\frac{4}{n}\right)\frac{1}{m}\sum_{j=1}^m\lam_j.\tag{HCY-2}
\end{equation*}
These two inequalities are known as Yang's first and second inequalities,
respectively. We note that \eqref{eq:HCY1} still holds if we replace
$\lam_{m+1}$ by an arbitrary $z\in(\lam_m,\lam_{m+1}]$ (see \cite{HaSt}), and
that the sharpest so far known explicit upper bound on $\lam_{m+1}$ is also
derived from \eqref{eq:HCY1}, see
\cite[formula (3.33)]{Ash}.

Payne-P\'olya-Weinberger, Hile-Protter and Yang inequalities are commonly referred
to as {\it universal estimates\/} for the eigenvalues of the Dirichlet Laplacian.
These estimates are valid uniformly over
all bounded domains in $\Rbb[n]$. The derivation of all four results is similar and
uses the variational principle with ingenious choices of test functions, and the
Cauchy-Schwarz inequality. We refer the reader to the extensive survey \cite{Ash}
which provides the detailed proofs as well as the
proof of the implication
$$
\text{\eqref{eq:HCY1}}\Longrightarrow\text{\eqref{eq:HCY2}}\Longrightarrow
\text{\eqref{eq:HP}}\Longrightarrow\text{\eqref{eq:PPW}}\,.
$$
In 1997, Harrell and Stubbe \cite{HaSt} showed that all of
these results are consequences of a certain abstract operator
identity and that this identity has several other applications.

Similar universal estimates were also obtained in spectral problems
for operators other then
the Euclidean Dirichlet Laplacian (or Schr\"odinger operator), e.g. higher order
differential operators in $\Rbb[n]$,
operators on manifolds, systems like Lam\'e system of elasticity etc., see,
\cite{Ha1, Ho1, Ho2, HaMi1, HaMi2} and already mentioned survey paper \cite{Ash}.

Unfortunately, despite the abstract nature of the results of \cite{HaSt},
it is unclear whether they are applicable in all these cases.

The first main result of our paper is a general abstract
operator identity which holds under minimal restrictions:

\begin{thm}\label{th:intro1}
Let $H$ and $G$ be self-adjoint operators such that
$G(D_H)\subseteq D_H$. Let $\lam_j$ and $\phi_j$ be eigenvalues
and eigenvectors of $H$. Then
for each $j$
\begin{equation}\label{eq:intro2}
\sum_k
\frac{|\scal{[H,G]\phi_j,\phi_k}|^2}{\lam_k-\lam_j}=
-\frac{1}{2}\scal{[[H,G],G]\phi_j,\phi_j}.
\end{equation}
\end{thm}

This theorem has a lot of applications, notably the estimates of the
eigenvalue gaps of various operators.
In particular, the results of Payne,
P\'olya and Weinberger for the Dirichlet Laplacian follow from 
\eqref{eq:intro2} if we
set $G$ to be an operator of multiplication by the coordinate $x_l$.
Then \eqref{eq:intro2} takes a particular simple and elegant form:
\begin{equation}\label{eq:intro3}
\sum_{k=1}^\infty \frac{\displaystyle\int_\Ome 
\frac{\dpar\phi_m}{\dpar x_l}\,\phi_k}{\lam_k-\lam_m}=\frac{1}{4}\,.
\end{equation}
(According to B~Simon \cite{Si}, this identity was known to physicists
already in the 1930s.)
Then \eqref{eq:PPW} follows from  \eqref{eq:intro3} if we sum the resulting equalities over $l$ and use some simple bounds, see
Examples~\ref{ex:vardirichlet},~\ref{ex:dirichlet} for details.
There are other
applications of Theorem \ref{th:intro1} -- in each particular case one
should work out what is the optimal choice of $G$ -- and we give
here several such applications.

Our other main main result is the generalization of the formula
\eqref{eq:intro2} to the case of several operators. Namely, suppose we
have two operators $H_1$ and $H_2$ (the model case being Laplacians
with different boundary conditions) and we want to estimate
eigenvalues of $H_1$ in terms of eigenvalues of $H_2$. Then one
can write the formula, similar to \eqref{eq:intro2}, but instead of
the usual commutator $[H,G]$ we will have the `mixing commutator'
$H_1G-GH_2$. It turns out that one of the operators $H_j$ in this
scheme can be non-self-adjoint. Details are given in Section
3. We give several applications of the second formula as well;
however, now the possible choice of the auxiliary operator $G$ is
even more restrictive, since we have to make sure that all the
commutators involved make sense.

\

{\bf Acknowledgements.\/} We are grateful to M Ashbauch, E B Davies,
and  E M Harrell for numerous helpful remarks and useful discussions.

\

\section{Statements for a Single Operator}

In this Section, $H$ denotes a self-adjoint operator with eigenvalues
$\lam_j$ and an orthonormal basis of eigenfunctions $\phi_j$.
Operator $H$ acts in a Hilbert space $\HC$ equipped with the scalar
product $\scal{\cdot,\cdot}$ and the corresponding norm
$\|\cdot\|$.

We start by stating the following obvious result.
\begin{lem}\label{lem:VL}
Let $\lam_j=\lam_k$. Then
\begin{equation}\label{eq:virial}
\scal{[H,G]\phi_j,\phi_k}=0\,.
\end{equation}
\end{lem}

Our next Theorem gives various trace identities similar to the one
given in Theorem~\ref{th:intro1}.

\begin{thm}\label{th:1}
Let $H$ and $G$ be self-adjoint operators with domains $D_H$ and $D_G$ such that
$G(D_H)\subseteq D_H\subseteq D_G$. Let $\lam_j$ and $\phi_j$ be eigenvalues
and eigenvectors of $H$. Let $P_j$ be the projector on the
eigenspace $\HC_j$ corresponding to the set of eigenvalues which
are equal to $\lam_j$. Then
for each $j$
\begin{equation}\label{eq:id1}
\sum_k
\frac{|\scal{[H,G]\phi_j,\phi_k}|^2}{\lam_k-\lam_j}=
-\frac{1}{2}\scal{[[H,G],G]\phi_j,\phi_j}.
\end{equation}
\begin{equation}\label{eq:id2}
\sum_k (\lam_k-\lam_j) |\scal{G\phi_j,\phi_k}|^2 =
-\frac{1}{2}\scal{[[H,G],G]\phi_j,\phi_j}.
\end{equation}
\begin{equation}\label{eq:id3}
\sum_k
\frac{|\scal{[H,G]\phi_j,\phi_k}|^2}{(\lam_k-\lam_j)^2} =
\|G\phi_j\|^2-\|P_jG\phi_j\|^2\,.
\end{equation}
\begin{equation}\label{eq:id4}
\sum_k (\lam_k-\lam_j)^2 |\scal{G\phi_j,\phi_k}|^2 = \|[H,G]\phi_j\|^2.
\end{equation}
\end{thm}

\begin{rem} The summation in (\ref{eq:id1})--(\ref{eq:id4}) is over all
$k$. Lemma~\ref{lem:VL} guarantees that the
summands in (\ref{eq:id1}) and (\ref{eq:id3}) are correctly
defined even when $\lam_k=\lam_j$ (if we assume $0/0=0$).
\end{rem}

\begin{rem}\label{rem:domains}
Instead of the condition $G(D(H))\subseteq D(H)$
we can impose weaker conditions $G\phi_j\in D(H)$, $G^2\phi_j\in D(H)$, $j=1,\dots$.
Moreover, the latter condition can be dropped if the double commutator is understood
in the weak sense, i.e., if the right-hand side of \eqref{eq:id1} and \eqref{eq:id2}
is replaced by $\scal{[H,G]\phi_j,G\phi_j}$ (see \eqref{eq:byparts2} below).
\end{rem}

\begin{rem}
Formulae (\ref{eq:id1})--(\ref{eq:id4}) can be extended to the case
of $H$ having continuous spectrum. In this case, the identities will include
integration instead of summation, cf. \cite{HaSt}. We omit the full details.
\end{rem}

\begin{proof}[Proof of Theorem~\ref{th:1}]

We are going to prove identities (\ref{eq:id1}) and (\ref{eq:id2});
the other two identities are proved in a similar manner (and are much easier).

Obviously, we have
\begin{equation}\label{eq:orthog}
[H,G]\phi_j=(H-\lam_j)G\phi_j\,.
\end{equation}

Therefore,

\begin{equation}\label{eq:old2}
\scal{G[H,G]\phi_j,\phi_j}=\scal{G(H-\lam_j)G\phi_j,\phi_j}\,.
\end{equation}

Since $G$ is self-adjoint, we have
\begin{equation}\label{eq:byparts1}
\begin{split}
&\scal{G(H-\lam_j)G\phi_j,\phi_j}=
\scal{(H-\lam_j)G\phi_j,G\phi_j}\\
&=\sum_k \scal{(H-\lam_j)G\phi_j,\phi_k}\scal{\phi_k,G\phi_j}
=\sum_k (\lam_k-\lam_j)|\scal{G\phi_j,\phi_k}|^2.
\end{split}
\end{equation}

Using the fact that $[H,G]$ is skew-adjoint, the left-hand side of
\eqref{eq:old2} can be rewritten as

\begin{equation}
\begin{split}
\scal{G[H,G]\phi_j,\phi_j}&=-\scal{[[H,G],G]\phi_j,\phi_j}+
\scal{[H,G]G\phi_j,\phi_j}\\&=-\scal{[[H,G],G]\phi_j,\phi_j}-
\scal{\phi_j,G[H,G]\phi_j}\,,
\end{split}
\end{equation}
so
\begin{equation}\label{eq:byparts2}
\scal{G[H,G]\phi_j,\phi_j}=-\frac{1}{2}\scal{[[H,G],G]\phi_j,\phi_j}
\end{equation}
(notice that $\scal{G[H,G]\phi_j,\phi_j}$ is
real, see \eqref{eq:old2} and \eqref{eq:byparts1}).
This proves \eqref{eq:id2}.

Since \eqref{eq:orthog} implies
$$
\scal{[H,G]\phi_j,\phi_k}=(\lam_k-\lam_j)\scal{G\phi_j,\phi_k}\,,
$$
this also proves \eqref{eq:id1}.
\end{proof}

Let us now put in \eqref{eq:id3} $G=[H,F]$ where $F$ is
skew-adjoint. Then due to
\eqref{eq:virial} the second term in the right-hand side of \eqref{eq:id3}
vanishes, and we have the following
\begin{cor} For a skew-adjoint operator $F$ such that $F(\phi_j)\in D(H^2)$ for all
$j$, we have
\begin{equation}\label{eq:id3bis}
\sum_k
\frac{|\scal{[H,[H,F]]\phi_j,\phi_k}|^2}{(\lam_k-\lam_j)^2} =
\|[H,F]\phi_j\|^2.
\end{equation}
\end{cor}

As above (see Remark~\ref{rem:domains}), we can replace the conditions
$F(\phi_j)\in D(H^2)$ by weaker ones $F(\phi_j)\in D(H)$ if we agree to
understand the double commutators in an appropriate weak sense.

From now on, we assume that the sequence of eigenvalues
$\{\lam_j\}_{j=1}^\infty$ is non-decreasing.

We now have at our disposal all the tools required for establishing
the ``abstract'' versions of \eqref{eq:PPW} and \eqref{eq:HCY1}.

\begin{cor}\label{cor:maincor}
Under conditions of Theorem~\ref{th:1},
\begin{equation}\label{eq:ourPPW}
-(\lam_{m+1}-\lam_m)\sum_{j=1}^m([[H,G],G]\phi_j,\phi_j)
\le 2\sum_{j=1}^m\|[H,G]\phi_j\|^2\,.
\end{equation}
\end{cor}

\begin{proof}
Let us sum the equations \eqref{eq:id1} over $j=1,...,m$. Then
we have
\begin{equation}\label{eq:3}
\sum_{j=1}^m\sum_{k=m+1}^\infty
\frac{|([H,G]\phi_j,\phi_k)|^2}{\lam_k-\lam_j}=
-\frac{1}{2}\sum_{j=1}^m([[H,G],G]\phi_j,\phi_j)\,.
\end{equation}

Parceval's
equality implies that the left-hand side of \eqref{eq:3} is not greater than
$\displaystyle\frac{1}{\lam_{m+1}-\lam_m}\sum_{j=1}^m\|[H,G]\phi_j\|^2$.
This proves \eqref{eq:ourPPW}.
\end{proof}

The next corollary uses the idea of \cite{HaSt}.

\begin{cor}\label{cor:maincor1}
For all $z\in(\lam_m,\lam_{m+1}]$ we have:
\begin{equation}\label{eq:cor1}
\sum_{j=1}^m
(z-\lam_j)\|[H,G]\phi_j\|^2\ge
-\frac{1}{2}\sum_{j=1}^m(z-\lam_j)^2\scal{[[H,G],G]\phi_j,\phi_j}
\end{equation}
\end{cor}

\begin{proof}
Let us multiply \eqref{eq:id1} by $(z-\lam_j)^2$ and sum the
result over all $j=1,...,m$. We will get:
\begin{equation}\label{eq:sumj}
\sum_{j=1}^m\sum_k
(z-\lam_j)^2\frac{|\scal{[H,G]\phi_j,\phi_k}|^2}{\lam_k-\lam_j}=
-\frac{1}{2}\sum_{j=1}^m(z-\lam_j)^2\scal{[[H,G],G]\phi_j,\phi_j}.
\end{equation}

The left-hand side of \eqref{eq:sumj} can be estimated as follows:

\begin{equation}\label{eq:sumj1}
\begin{split}
&\sum_{j=1}^m\sum_k
(z-\lam_j)^2\frac{|\scal{[H,G]\phi_j,\phi_k}|^2}{\lam_k-\lam_j}\\
&\quad=\sum_{j=1}^m\sum_{k=1}^m
(z-\lam_j)^2\frac{|\scal{[H,G]\phi_j,\phi_k}|^2}{\lam_k-\lam_j}\\
&\qquad+\sum_{j=1}^m\sum_{k=m+1}^\infty (z-\lam_j)^2\frac{|\scal{[H,G]\phi_j,\phi_k}|^2}{\lam_k-\lam_j}\\
&\quad\le \sum_{j=1}^m\sum_{k=1}^m
(z-\lam_j)^2\frac{|\scal{[H,G]\phi_j,\phi_k}|^2}{\lam_k-\lam_j}\\
&\qquad+\sum_{j=1}^m\sum_{k=m+1}^\infty (z-\lam_j)|\scal{[H,G]\phi_j,\phi_k}|^2\\
&\quad=\sum_{j=1}^m (z-\lam_j)\sum_{k=1}^\infty|\scal{[H,G]\phi_j,\phi_k}|^2\\
&\quad+\sum_{j=1}^m\sum_{k=1}^m\left(
(z-\lam_j)^2\frac{|\scal{[H,G]\phi_j,\phi_k}|^2}{\lam_k-\lam_j}
-(z-\lam_j)|\scal{[H,G]\phi_j,\phi_k}|^2\right)\\
&\quad=\sum_{j=1}^m
(z-\lam_j)\|[H,G]\phi_j\|^2\\
&\qquad+\sum_{j=1}^m\sum_{k=1}^m\left(
(z-\lam_j)|\scal{[H,G]\phi_j,\phi_k}|^2\left(\frac{z-\lam_j}{\lam_k-\lam_j}
-1\right)\right)\\
&\quad=
\sum_{j=1}^m
(z-\lam_j)\|[H,G]\phi_j\|^2\\
&\qquad+\sum_{j=1}^m\sum_{k=1}^m\left(
\frac{(z-\lam_j)(z-\lam_k)}{\lam_k-\lam_j}|\scal{[H,G]\phi_j,\phi_k}|^2\right)\\
&\quad=
\sum_{j=1}^m
(z-\lam_j)\|[H,G]\phi_j\|^2\,.
\end{split}
\end{equation}
(The last equality uses the fact that the expression under
$\displaystyle \sum_{j=1}^m\sum_{k=1}^m$ is skew-symmetric with respect to
$j,k$.) Now \eqref{eq:sumj} and \eqref{eq:sumj1} imply
\eqref{eq:cor1}.
\end{proof}

\begin{rem} As we will see in case of the Dirichlet Laplacian, our formula
\eqref{eq:ourPPW} is an abstract generalization of Payne-P\'olya-Weinberger
formula \eqref{eq:PPW}, and \eqref{eq:cor1} is an abstract generalization of
Yang's formula \eqref{eq:HCY1}.
\end{rem}

\section{Statements for a Pair of Operators}

The results of previous Section are not applicable, directly, to
non-self-adjoint operators. To extend
the spectral trace identities to a non-self-adjoint case we consider
{\bfseries pairs of operators}
$H_1$, $H_2$, where one of them is allowed to be non-self-adjoint.
Using auxiliary operators $G_1$, $G_2$, we can relate the spectra of
$H_1$ and $H_2$.

First, we introduce the following notation. For a {\bf triple of operators}
$X$, $Y$, $Z$ acting in a Hilbert space $\HC$ we define the ``mixing commutators''
\begin{equation}\label{eq:triplecomm}
[X,Y;Z]=XZ-ZY\,,\qquad\{X,Y;Z\}_\pm=XZ\pm Z^*Y\,.
\end{equation}

We note some elementary properties of ``mixing commutators'' \eqref{eq:triplecomm}:
\begin{gather*}
[X,X;Z]=[X,Z]\,,\qquad [X,Y;Z]^*=-[Y^*,X^*;Z^*]\,,\\
\{X,Y;Z\}_\pm^*=\pm\{Y^*,X^*;Z\}_\pm\,.
\end{gather*}

We always assume non-self-adjoint operators
to be closed.

Our main result concerning non-self-adjoint operators is the following

\begin{thm}\label{th:nonsa1}
Let $H_1$ be a self-adjoint operator in a Hilbert space $\HC$ with eigenvalues
$\lam_k$ and an orthonormal basis of eigenfunctions $\phi_k$, and let $H_2$ be a
(not necessarily self-adjoint) operator in $\HC$ with eigenvalues $\mu_j$
and eigenfunctions $\psi_j$. Define, for an auxiliary pair of
operators $G_1$, $G_2$ in $\HC$, the operators
\begin{equation}\label{eq:triplecommABCD}
\begin{split}
A&=[H_1,H_2;G_1^*]\,,\\
B&=[H_1,H_2;G_2]\,,\\
C&=[H_2^*,H_1;G_2^*]=-B^*\,,\\
D_\pm&={\{C,B;G_1^*\}}_\pm\,.
\end{split}
\end{equation}

If the operators $A$, $B$, and $D_\pm$ are well defined, and all the eigenfunctions
of $H_2$ belong to their domains, then the
following trace identities hold for any fixed $j$:
\begin{align}
\re\sum_k\frac{\lam_k-\mu_j}{|\lam_k-\mu_j|^2}\scal{B\psi_j,\phi_k}\cdot
\overline{\scal{A\psi_j,\phi_k}}&=-\frac{1}{2}\scal{D_{-}\psi_j,\psi_j}\,,
\label{eq:nonsati1}
\\
i\im \sum_k\frac{\lam_k-\mu_j}{|\lam_k-\mu_j|^2}\scal{B\psi_j,\phi_k}\cdot
\overline{\scal{A\psi_j,\phi_k}}&=\frac{1}{2}\scal{D_{+}\psi_j,\psi_j}\,.
\label{eq:nonsati2}
\end{align}
\end{thm}

\begin{proof} Acting as in the proof of Theorem~\ref{th:1} we get
\begin{equation}\label{eq:nonsapr1}
\begin{split}
\scal{G_1[H_1,H_2;G_2]\psi_j,\psi_j}&=\scal{G_1(H_1G_2-G_2H_2)\psi_j,\psi_j}\\
&=\scal{(H_1-\mu_j)G_2\psi_j,G_1^*\psi_j}\\
&=\sum_k\scal{(H_1-\mu_j)G_2\psi_j,\phi_k}\cdot\scal{\phi_k,G_1^*\psi_j}\\
&=\sum_k\scal{G_2\psi_j,(H_1-\overline{\mu_j})\phi_k}\cdot\scal{\phi_k,G_1^*\psi_j}\\
&=\sum_k(\lam_k-\mu_j)\overline{\scal{G_1^*\psi_j,\phi_k}}\cdot\scal{G_2\psi_j,\phi_k}\,.
\end{split}
\end{equation}

Also,
\begin{equation}\label{eq:nonsapr2}
\begin{split}
\scal{[H_1,H_2;G_2]\psi_j,\phi_k}&=\scal{(H_1G_2-G_2H_2)\psi_j,\phi_k}\\
&=\lam_k\scal{G_2\psi_j,\phi_k}-\scal{G_2\mu_j\psi_j,\phi_k}\\
&=(\lam_k-\mu_j)\scal{G_2\psi_j,\phi_k}\,,
\end{split}
\end{equation}
and, similarly,
\begin{equation}\label{eq:nonsapr3}
\overline{\scal{[H_1,H_2;G_1^*]\psi_j,\phi_k}}=
(\lam_k-\overline{\mu_j})\overline{\scal{G_1^*\psi_j,\phi_k}}\,.
\end{equation}
Therefore, \eqref{eq:nonsapr1} can be re-written as
\begin{multline}\label{eq:nonsapr4}
\scal{G_1[H_1,H_2;G_2]\psi_j,\psi_j}\\=
\sum_k\frac{\lam_k-\mu_j}{|\lam_k-\mu_j|^2}
\scal{[H_1,H_2;G_2]\psi_j,\phi_k}\cdot
\overline{\scal{[H_1,H_2;G_1^*])\psi_j,\phi_k}}\,.
\end{multline}

Finally, using the definitions \eqref{eq:triplecomm}, we have
\begin{equation}\label{eq:nonsapr5}
\begin{split}
2\re\scal{G_1[H_1,H_2;G_2]\psi_j,\psi_j}&=
\scal{(G_1[H_1,H_2;G_2]+[H_1,H_2;G_2]^*G_1^*)\psi_j,\psi_j}\\
&=-\scal{(-G_1[H_1,H_2;G_2]+[H_2^*,H_1;G_2^*])G_1^*\psi_j,\psi_j}\\
&=-\scal{\{[H_2^*,H_1;G_2^*],[H_1,H_2;G_2];G_1^*\}_-\psi_j,\psi_j}\,.
\end{split}
\end{equation}
and
\begin{equation}\label{eq:nonsapr6}
\begin{split}
2i\im\scal{G_1[H_1,H_2;G_2]\psi_j,\psi_j}&=
\scal{(G_1[H_1,H_2;G_2]-[H_1,H_2;G_2]^*G_1^*)\psi_j,\psi_j}\\
&=\scal{(G_1[H_1,H_2;G_2]+[H_2^*,H_1;G_2^*]G_1^*)\psi_j,\psi_j}\\
&=\scal{\{[H_2^*,H_1;G_2^*],[H_1,H_2;G_2];G_1^*\}_+\psi_j,\psi_j}\,.
\end{split}
\end{equation}

The Theorem now follows by combining \eqref{eq:nonsapr4}-\eqref{eq:nonsapr6}
and using \eqref{eq:triplecommABCD}.
\end{proof}

The trace identities \eqref{eq:nonsati1}, \eqref{eq:nonsati2} are much
simpler if we choose $G_2^*=G_1$. Then $A=B=[H_1,H_2;G_1^*]$,
and we immediately arrive at

\begin{thm}\label{th:nonsa2} If, in addition to conditions of
Theorem~\ref{th:nonsa1}, we assume $G_2^*=G_1$, the following trace
identities hold for any $j$,
\begin{align}
\sum_k\frac{\lam_k-\re\mu_j}{|\lam_k-\mu_j|^2}|\scal{A\psi_j,\phi_k}|^2&=
-\frac{1}{2}\scal{{\{-A^*,A;G_1^*\}}_{-}\psi_j,\psi_j}\,,
\label{eq:nonsati3}
\\
i\sum_k\frac{\im\mu_j}{|\lam_k-\mu_j|^2}|\scal{A\psi_j,\phi_k}|^2&=
\frac{1}{2}\scal{{\{-A^*,A;G_1^*\}}_{+}\psi_j,\psi_j}\,.
\label{eq:nonsati4}
\end{align}
\end{thm}

An even simpler case is when the operators $H_2$ and $G_1=G_2$ are self-adjoint.
As for any self-adjoint $Z$, $\{X,Y;Z\}_-=[X,Y;Z]$, we do not have to use any
``curly brackets'' commutators and immediately obtain

\begin{thm}\label{th:nonsa3} If, in addition to conditions of
Theorem~\ref{th:nonsa1}, we assume that $H_2=H_2^*$ and $G_1=G_1^*=G_2=G$,
the following
trace identity holds for any $j$:
\begin{equation}\label{eq:nonsati5}
\sum_k\frac{1}{\lam_k-\mu_j}|\scal{[H_1,H_2;G]\psi_j,\phi_k}|^2=
-\frac{1}{2}\scal{[[H_2,H_1;G],[H_1,H_2;G];G]\psi_j,\psi_j}\,.
\end{equation}
\end{thm}

We emphasize that each of the Theorems~\ref{th:nonsa1}--\ref{th:nonsa3} supersedes
Theorem~\ref{th:1}. Indeed, if we set $H_1=H_2=H$, $\mu_k=\lambda_k$,
$\psi_k=\phi_k$, and
$G_1=G_2=G$, we have $[H,H;G]=[H,G]$, $[[H,H;G], [H,H;G]; G]=[[H,G],G]$, and
identity
\eqref{eq:nonsati5} becomes \eqref{eq:id1}. The other identities generalizing
\eqref{eq:id2}--\eqref{eq:id4} in Theorem~\ref{th:1},
can be obtained in similar fashion.

\begin{rem} The main difficulty in applying Theorems~\ref{th:nonsa1}--\ref{th:nonsa3}
is the choice of auxiliary operators $G_1$ and $G_2$ in such a way that all
the commutators involved make sense.
Similarly to Remark~\ref{rem:domains}, we can weaken the conditions of the
Theorems by considering the double ``mixing'' commutators in the weak sense only.
\end{rem}

In principle, one can obtain estimates for the eigenvalues in a
general situation of Theorem~\ref{th:nonsa1}. However, this is impractical because
of the variety of combinations of signs of terms in
\eqref{eq:nonsati1} and \eqref{eq:nonsati2}. The situation simplifies
if we consider more restricted choice of Theorems~\ref{th:nonsa2} and \ref{th:nonsa3}.

We start with applications of Theorem~\ref{th:nonsa2}. Before stating the main results
we introduce the following notation in addition to \eqref{eq:triplecommABCD}:
\begin{equation}\label{eq:morenotation}
a_j=\|A\psi_j\|^2\,,\qquad d_j^-=-\scal{D_-\psi_j,\psi_j},,\qquad
d_j^+=-i\scal{D_+\psi_j,\psi_j}
\end{equation}
(recall that $A=[H_1,H_2;G_1^*]$, $D_\pm=\{-A^*,A;G_1^*\}_\pm$).
It is easy to check that $d_j^\pm$ are in fact real numbers.

\begin{cor}\label{cor:nonsaest1} Under conditions of Theorem~\ref{th:nonsa2},
for any fixed $j$,
\begin{equation}\label{eq:nonsaest1}
\dist(\mu_j,\spec H_1)\le\frac{2a_j}{\sqrt{(d_j^-)^2+(d_j^+)^2}}\,.
\end{equation}
Moreover,
\begin{equation}\label{eq:nonsaest2}
\min_k |\re\mu_j-\lambda_k|\le\min_k\frac{|\mu_j-\lambda_k|^2}{|\re\mu_j-\lambda_k|}
\le\frac{2a_j}{|d_j^-|}
\end{equation}
and
\begin{equation}\label{eq:nonsaest3}
|\im\mu_j|\le\min_k\frac{|\mu_j-\lambda_k|^2}{|\im\mu_j|}
\le\frac{2a_j}{|d_j^+|}\,.
\end{equation}
\end{cor}

\begin{proof} Subtracting identity \eqref{eq:nonsati3} from
\eqref{eq:nonsati4}, taking the absolute value, and using the triangle
inequality and \eqref{eq:morenotation}, we get
$$
\sum_k \frac{1}{|\lambda_k-\mu_j|}|\scal{A\psi_j,\phi_k}|^2
\ge\frac{1}{2}|d_j^-+id_j^+|\,.
$$
The left-hand side of this inequality is estimated from above by
\begin{equation*}
\begin{split}
\max_k\frac{1}{|\lam_k-\mu_j|}\sum_k|\scal{A\psi_j,\phi_k}|^2&=
\frac{1}{\displaystyle \min_k |\mu_j-\lambda_k|}\|A\psi_j\|^2\\
&=\frac{1}{\dist(\mu_j,\spec H_1)}a_j\,,
\end{split}
\end{equation*}
which implies \eqref{eq:nonsaest1}. The estimates \eqref{eq:nonsaest2} and
\eqref{eq:nonsaest3}
are obtained by applying exactly the same procedure to \eqref{eq:nonsati3}
and \eqref{eq:nonsati4}
separately.
\end{proof}

\section{Examples}
\begin{exmp}\label{ex:vardirichlet}
{\bfseries Second order operator with variable coefficients, Dir\-ich\-let
problem\/}.
Let $\dpar_k=\dpar/\dpar x_k$, and let
$H=-\sum_{k,l=1}^n \dpar_k a_{kl}(x)\dpar_l$ be a positive elliptic operator
with Dirichlet boundary conditions in $\Ome\subset\Rbb[n]$ ($A=\{a_{jk}\}$
is positive).
Let $G$ be an operator of multiplication by a function $f$. Then
$$
[H,G]u=(Hf)u-2\sum_{k,l=1}^n (\dpar_k f) a_{kl}(x)(\dpar_l u)\,,
$$
and
$$
[[H,G],G]=-2\sum_{k,l=1}^n (\dpar_k f) a_{kl}(x)(\dpar_l f)\,.
$$

Therefore, Corollary \ref{cor:maincor} implies:

\begin{equation}\label{eq:vari1}
\lam_{m+1}-\lam_m
\le\frac{\displaystyle\sum_{j=1}^m\int_\Om\left((Hf)\phi_j-
2\sum_{k,l=1}^n (\dpar_k f) a_{kl}(x)(\dpar_l \phi_j)\right)^2}
{\displaystyle\sum_{j=1}^m \int_\Om\sum_{k,l=1}^n(\dpar_k f) a_{kl}(x)(\dpar_l
f)\phi_j^2}
\end{equation}

Now, each choice of $f$ in \eqref{eq:vari1} will produce an
inequality for the spectral gap. For example, we can choose
$f=x_i$. Then \eqref{eq:vari1} will have the following form:

\begin{equation}\label{eq:vari2}
\displaystyle
\lam_{m+1}-\lam_m
\le \frac{\displaystyle\sum_{j=1}^m\int_\Om\left(\sum_{l=1}^n(\dpar_{l}a_{li}(x))\phi_j+
2\sum_{l=1}^n  a_{il}(x)(\dpar_l \phi_j)\right)^2}
{\displaystyle\int_\Om a_{ii}(x) \sum_{j=1}^m\phi_j^2}\,.
\end{equation}

Since \eqref{eq:vari2} is valid for all $i$, we have:

\begin{equation}\label{eq:vari3}
\begin{split}
\lam_{m+1}-\lam_m
&\le \frac{\displaystyle\sum_{i=1}^n\sum_{j=1}^m
\int_\Om\left(\sum_{l=1}^n(\dpar_{l}a_{li}(x))\phi_j+
2\sum_{l=1}^n  a_{il}(x)(\dpar_l \phi_j)\right)^2}
{\displaystyle\sum_{j=1}^m \int_\Om \Tr(A(x))\phi_j^2}\\
&\le\frac{\displaystyle p\sum_{i=1}^n\sum_{j=1}^m
\int_\Om\left(\sum_{l=1}^n(\dpar_{l}a_{li}(x))\right)^2\phi_j^2}
{\displaystyle\sum_{j=1}^m \int_\Om \Tr(A(x))\phi_j^2}\\
&+
\frac{\displaystyle 4q\sum_{i=1}^n\sum_{j=1}^m
\int_\Om\left(\sum_{l=1}^n  a_{il}(x)(\dpar_l \phi_j)\right)^2}
{\displaystyle\sum_{j=1}^m \int_\Om \Tr(A(x))\phi_j^2}\,,
\end{split}
\end{equation}
where $p$ and $q$ are arbitrary positive numbers greater than one
such that $(p-1)(q-1)=1$.
The first term in the right-hand side of \eqref{eq:vari3} can be estimated by

\begin{equation}
\sup_{x\in\Om}
\frac{\displaystyle p\sum_{i=1}^n\left(\sum_{l=1}^n(\dpar_{l}a_{li}(x))\right)^2}
{\displaystyle m\Tr(A(x))}\,.
\end{equation}
The second term is not greater than

\begin{equation}
\frac{\displaystyle 4q(\sum_{j=1}^m \lam_j)
\sup_{x\in\Om}(\text{maximal eigenvalue of }A(x))}
{\displaystyle m\inf_{x\in\Om}\Tr(A(x))}.
\end{equation}

This gives the inequality for the spectral gap:

\begin{equation}\label{eq:pq}
\begin{split}
\lam_{m+1}-\lam_m
&\le
\sup_{x\in\Om}
\frac{\displaystyle p\sum_{i=1}^n\left(\sum_{l=1}^n(\dpar_{l}a_{li}(x))\right)^2}
{m\Tr(A(x))}\\
&+
\frac{\displaystyle 4q\left(\sum_{j=1}^m \lam_j\right)
\sup_{x\in\Om}(\text{maximal eigenvalue of }A(x))}
{\displaystyle m\inf_{x\in\Om}\Tr(A(x))}
\end{split}
\end{equation}
in terms of the previous eigenvalues and properties of the
coefficients of the operator but not the geometric characteristics
of the domain.
\end{exmp}

\begin{exmp}\label{ex:dirichlet}
{\bfseries Dirichlet Laplacian\/}. Let now $H=-\Del$ acting in the bounded domain
$\Om\subset\Rbb[n]$ with Dirichlet boundary conditions. Then
in \eqref{eq:pq} we can let $p\to\infty$ (and so $q\to 1$) and get
\eqref{eq:PPW} inequality (in the same way as in \cite{HaSt}):

\begin{equation}
\lam_{m+1}-\lam_m\le \frac{4}{mn}\sum_{j=1}^m \lam_j.
\end{equation}

If one uses Corollary~\ref{cor:maincor1} instead, one gets the
following inequality (in the same way as in \cite{HaSt})
for all $z\in(\lam_m,\lam_{m+1}]$:

\begin{equation}\label{eq:cor1eq1}
\frac{4}{n}\sum_{j=1}^m
(z-\lam_j)\lam_j\ge
\sum_{j=1}^m(z-\lam_j)^2.
\end{equation}

If $z=\lam_{m+1}$, \eqref{eq:cor1eq1} becomes \eqref{eq:HCY1}.

Now let us look once again at our main identity when $H$ is the
Dirichlet Laplacian and $G$ is the operator of multiplication by $x_l$
($l=1,...,n$):

\begin{equation}\label{eq:dirichlet}
\sum_{k=1}^\infty \frac{w_{m,k,l}^2}{\lam_k-\lam_m}=\frac{1}{4},
\end{equation}
where
\begin{equation}
w_{m,k,l}:=\int_{\Om}\frac{\partial\phi_m}{\partial x_l}\phi_k.
\end{equation}
Using Gaussian elimination, one can find the orthogonal
coordinate system $x_1,...,x_n$ such that
\begin{equation}
\begin{split}
w_{m,m+1,1}&=w_{m,m+1,2}=...=w_{m,m+1,n-1}=w_{m,m+2,1}=...\\
&=w_{m,m+2,n-2}=...=w_{m,m+n-1,1}=0.
\end{split}
\end{equation}
Let us now make the obvious
estimate of the left-hand side of \eqref{eq:dirichlet}:
\begin{equation}
\frac{1}{\lam_{m+l}-\lam_m}
\int_\Om \biggl(\frac{\partial\phi_m}{\partial x_l}\biggr)^2\ge
\sum_{k=1}^\infty \frac{w_{m,k,l}^2}{\lam_k-\lam_m}=\frac{1}{4},
\end{equation}
or
\begin{equation}
\lam_{m+l}-\lam_m\le 4\int_\Om \biggl(\frac{\partial\phi_m}{\partial
x_l}\biggr)^2.
\end{equation}
Summing these inequalities over all $l=1,...,n$ gives
\begin{equation}
\sum_{l=1}^n\lam_{m+l}\le (4+n)\lam_m.
\end{equation}
As far as we know, this estimate is new for $m>1$ (for a discussion of
the case $m=1$ see \cite[Section 3.2]{Ash}).
\end{exmp}

\begin{exmp}
{\bfseries Neumann Laplacian\/}.
The case of the Neumann conditions is much more difficult than the
Dirichlet ones because now if we take $G$ to be a multiplication
by a function $g$, we have to make sure that $g$ satisfies Neumann
conditions on the boundary. Therefore, we cannot get any
eigenvalue estimates without the preliminary knowledge of the
geometry of $\Om\subset\Rbb[n]$. We combine the ideas of
\cite{HaMi1} and \cite{ChGrYa} to get some improvement on the estimate of
\cite{HaMi1}.

Suppose, for example, that that we can
insert $q$ balls $B_p=B(x_p,r_p)$ ($p=1,...,q$)
of radii $r_1\ge r_2\ge\dots\ge r_q$ inside $\Om$
such that these balls do not intersect each other. Let $R(x)$
be the second radial eigenfunction of the Neumann Laplacian in a unit
ball $B(0,1)$ normalized in such a way that it is equal to $1$ on the
boundary of the ball. Then the function
\begin{equation}\label{eq:g}
g(x):=\begin{cases}
R(r_p^{-1}(x-x_p)),& x\in B_p\\
1,&\text{otherwise}
\end{cases}
\end{equation}
satisfies Neumann conditions on $\partial\Om$. Therefore, if we
take $G$ to be multiplication by $g$ and $H$ to be Neumann
Laplacian on $\Om$, they satisfy conditions of \ref{th:1}.
Therefore, corollary \ref{cor:maincor} implies (by $C_1,C_2,...$ we denote
different constants depending only on $n$)
\begin{multline}\label{eq:g1}
\lam_{m+1}-\lam_m\\
\le \frac{\displaystyle C_1\sum_{j=1}^m\sum_{p=1}^q
r_p^{-4}\int_{B_p}\phi_j^2R_p^2+C_2
\sum_{j=1}^m\sum_{p=1}^q\int_{B_p}|\nabla\phi_j|^2|\nabla R_p|^2}
{\displaystyle \sum_{j=1}^m\sum_{p=1}^q\int_{B_p}\phi_j^2|\nabla R_p|^2}.
\end{multline}
The denominator in the right-hand side of \eqref{eq:g1} can be estimated from the
below by noticing that $\phi_1\equiv \frac{1}{|\Om|}$. Therefore,
\begin{equation}\label{eq:g2}
\lam_{m+1}-\lam_m\le \frac{C_3|\Om|}{\displaystyle \sum_{p=1}^q
r_p^{-2+n}}\left(\sum_{p=1}^q r_p^{-4}+
r_q^{-2}\sum_{j=1}^m\lam_j\right).
\end{equation}
Assuming that all the radii $r_j$ are the same, we get
\begin{equation}\label{eq:g3}
\lam_{m+1}-\lam_m\le C_4|\Om|
r_q^{-n}\left(r_q^{-2}+
\frac{1}{q}\sum_{j=1}^m\lam_j\right).
\end{equation}
\end{exmp}

\begin{exmp}
{\bfseries Elasticity\/}. Here we mostly follow the lines of \cite{Ho2}
(though the final result is slightly different); for convenience we use
the same notation. We consider the spectral problem for the operator
of linear elasticity,
\begin{equation}\label{eq:elast0}
H\ubold=-\Del\ubold-\alp\bgrad\diver\ubold
\end{equation}
on a compact domain $\Ome\subset\Rbb[n]$ with smooth boundary, with Dirichlet
boundary conditions $\ubold|_{\partial\Ome}=\mathbf{0}$. Here
$\ubold=(u_1,\dots,u_n)$ is an $n$-dimensional vector-function of
$\xbold=(x_1,\dots,x_n)\in\Ome$, and $\alpha>0$
is a fixed parameter. Denote the eigenvalues of \eqref{eq:elast0}
by $\Lam_1\le\Lam_2\le\dots\Lam_j\le\dots$, and corresponding
eigenvectors $\ubold_j$.

We denote $L=-\Delta$, $M=-\bgrad\diver$, so that $H=L+\alp M$, and consider
the operators $G_l$ of multiplication by $x_l$, $l=1,\dots,n$. Then, by
\cite[Lemmas~4,~5]{Ho2},  we have
$$
[L, G_l]=-2S_l\,,\qquad [M,G_l]=-R_l\,,
$$
where
$S_l\ubold=\frac{\dpar\ubold}{\dpar x_l}$,
$R_l\ubold=(\diver\ubold)\bgrad x_l+\bgrad u_l$.
Also,
$$
\sum_{l=1}^n [R_l,G_l]\ubold=2\ubold\,,\qquad
\sum_{l=1}^n [S_l,G_l]\ubold=n\ubold\,.
$$

Applying the identity (\ref{eq:id1}) of Theorem~\ref{th:1} with $G=G_l$ and
summing over $l=1\dots n$, we obtain
$$
\sum_k\frac{\sum_{l=1}^n |\langle(2S_l+\alp R_l)\ubold_j,\ubold_k\rangle|^2}%
{\Lam_k-\Lam_j}=(n+\alp)\,.
$$

Corollary~\ref{cor:maincor} now implies the estimate
\begin{equation}\label{eq:elast1}
\Lam_{m+1}-\Lam_m \le
\frac{1}{m(n+\alp)}\sum_{j=1}^m\sum_{l=1}^n\|(2S_l+\alp R_l)\ubold_j\|^2\,.
\end{equation}
To estimate the right-hand side of (\ref{eq:elast1}), we need the following

\begin{lem}\label{lem:elast}
If $\ubold=\mathbf{0}$ on $\dpar\Ome$, then
\begin{gather}
\langle-\bgrad\diver\ubold,\ubold\rangle=\|\diver\ubold\|^2\,,
\label{eq:elast2}\\
\sum_{l=1}^n\|R_l\ubold\|^2=(n+2)\langle-\bgrad\diver\ubold,\ubold\rangle+
\langle-\Del\ubold,\ubold\rangle\,,
\label{eq:elast3}\\
\sum_{l=1}^n\|S_l\ubold\|^2=\langle-\Del\ubold,\ubold\rangle\,,
\label{eq:elast4}\\
\sum_{l=1}^n\langle S_l\ubold, R_l\ubold\rangle=
2\langle-\bgrad\diver\ubold,\ubold\rangle\,.
\label{eq:elast5}
\end{gather}
\end{lem}

\begin{proof}[Proof of Lemma~\ref{lem:elast}]
The equalities (\ref{eq:elast2})--(\ref{eq:elast4}) are proved in \cite{Ho2};
it remains only to prove (\ref{eq:elast5}).

Using the definitions of $R_l$, $S_l$, and integrating by parts, we have
\begin{equation*}
\begin{split}
\sum_{l=1}^n \langle S_l\ubold, R_l\ubold\rangle &=
\sum_{l=1}^n \int_\Ome \left(\frac{\dpar\ubold}{\dpar x_l}\right)\cdot
\left((\diver\ubold)\bgrad x_l+\bgrad u_l\right)\\
&=\sum_{l=1}^n \int_\Ome (\diver\ubold)\frac{\dpar u_l}{\dpar x_l}+
\sum_{l=1}^n\sum_{k=1}^n\int_\Ome \frac{\dpar u_l}{\dpar x_k} \frac{\dpar u_k}{\dpar x_l}
\\
&=\int_\Ome (\diver \ubold)^2-\int_\Ome (\ubold\cdot\bgrad\diver\ubold)\\
&=-2\langle\bgrad\diver\ubold,\ubold\rangle\,.
\end{split}
\end{equation*}
\end{proof}

Applying now Lemma~\ref{lem:elast} to the right-hand side of \eqref{eq:elast1},
we have
\begin{equation*}
\begin{split}
\Lam_{m+1}-\Lam_m &\le \frac{1}{m(n+\alp)}\sum_{j=1}^m
\left(4\|S_l\ubold_j\|^2+\alp^2\|R_l\ubold_j\|^2+4\alp\langle S_l\ubold_j, R_l\ubold_j\rangle\right)\\
&= \frac{1}{m(n+\alp)}\sum_{j=1}^m
\left((4+\alp^2)\langle-\Del\ubold_j,\ubold_j\rangle\right.\\
&\qquad\qquad\left.+((n+2)\alp^2+8\alp)\langle-\bgrad\diver\ubold_j,\ubold_j\rangle\right)\\
\le \frac{1}{m(n+\alp)}&\sum_{j=1}^m \max(4+\alp^2,(n+2)\alp+8)
\langle-\Del\ubold_j-\alp\bgrad\diver\ubold_j,\ubold_j\rangle\\
&=\frac{1}{m(n+\alp)}\max(4+\alp^2,(n+2)\alp+8)\sum_{j=1}^m \Lam_j\,.
\end{split}
\end{equation*}
\end{exmp}

\begin{exmp}
{\bfseries Two Schr\"odinger operators\/}. Here we consider a simple
example illustrating the results on pairs of operators. Let $H_1$ be
a Schr\"odinger operator $\displaystyle -\frac{d^2}{dx^2}+V_1(x)$ with
Neumann boundary conditions on a finite interval $I\subset\Rbb$ and $H_2$
be a Schr\"odinger operator $\displaystyle -\frac{d^2}{dx^2}+V_2(x)$ with
Dirichlet boundary conditions on the same interval; we assume that both
potentials are sufficiently smooth and that $V_1$ (but not necessarily
$V_2$) is real-valued.

We choose $G=G^*=G_1=G_2=i\frac{d}{dx}$. It easy to check that for
an eigenfunction $\psi$ of $H_2$ corresponding to an eigenvalue $\mu$ we have
$$
\left.(\frac{d}{dx})G\psi\right|_{\partial I}=
\left.i\frac{d^2}{dx^2}\psi\right|_{\partial I}=
-i (\mu-V_2)\psi|_{\partial I}=0\,.
$$
Thus, $G\psi\in D_{H_1}$, and the commutators appearing in
Theorem~\ref{th:nonsa2} are correctly defined.

Elementary computations then produce
$$
A=[H_1, H_2, G]=(V_1-V_2)i\frac{d}{dx}-iV_2'\,,\qquad
A^*=(V_1-\overline{V_2})+iV_1'\,,
$$
and, further on,
\begin{align}
D_+&=-A^*G+GA=(2i\im V_2)\frac{d^2}{dx^2}+2V_2'\frac{d}{dx}+V_2''\,,\label{eq:dplus}\\
D_-&=-A^*G-GA=2(V_1-\re V_2)\frac{d^2}{dx^2}+2(V_1'-V_2')\frac{d}{dx}-V_2''\,.\label{eq:dminus}
\end{align}

Substituting this expressions into \eqref{eq:nonsati3} and  \eqref{eq:nonsati4}, we obtain the trace identities,
\begin{multline}
\sum_k\frac{\lam_k-\re\mu_j}{|\lam_k-\mu_j|^2}|\scal{((V_1-V_2)i\frac{d}{dx}-iV_2')\psi_j,\phi_k}|^2\\
\quad=-\frac{1}{2}\scal{((2i\im V_2)\frac{d^2}{dx^2}+2V_2'\frac{d}{dx}+V_2'')\psi_j,\psi_j}\,,
\label{eq:shroid1}
\end{multline}
\begin{multline}
i\sum_k\frac{\im\mu_j}{|\lam_k-\mu_j|^2}|\scal{((V_1-V_2)i\frac{d}{dx}-iV_2')\psi_j,\phi_k}|^2\\
\quad=
\frac{1}{2}\scal{(2(V_1-\re V_2)\frac{d^2}{dx^2}+2(V_1'-V_2')\frac{d}{dx}-V_2'')\psi_j,\psi_j}\,.
\label{eq:shroid2}
\end{multline}

Also, the estimates \eqref{eq:nonsaest1}--\eqref{eq:nonsaest3} hold.

As usual, obtaining ``practical'' information about eigenvalues and eigenvalue gaps from
\eqref{eq:nonsaest1}--\eqref{eq:nonsaest3} requires constructing effective estimates
from above for
$$
a_j=\|A\psi_j\|^2=\|((V_1-V_2)i\frac{d}{dx}-iV_2')\psi_j\|^2\,,
$$
and from below for
$$
d_j^+=-i\scal{D_+\psi_j,\psi_j}=
-i\scal{((2i\im V_2)\frac{d^2}{dx^2}+2V_2'\frac{d}{dx}+V_2'')\psi_j,\psi_j}
$$
and
$$
d_j^-=-\scal{D_-\psi_j,\psi_j}=
-\scal{(2(V_1-\re V_2)\frac{d^2}{dx^2}+2(V_1'-V_2')\frac{d}{dx}-V_2'')\psi_j,\psi_j}\,.
$$

Estimating $a_j$ is easy:
$$
|a_j|\le \|V_1-V_2\|_1^2\lambda_j^2+\|V_2'\|_1^2\,,
$$
where $\|\cdot\|_1$ stands for the $L_1$ norm on the interval.

The estimation of $d_j^\pm$ doesn't seem to be possible in general, without
additional assumptions on
potentials $V_1$ and $V_2$. Therefore, we shall consider a simple particular
case of $V_1=V_2=V$, assuming additionally that $V''\ge c>0$ uniformly on $I$.
Then we have
$$
a_j=\|V'\psi_j\|^2\le\|V'\|_1^2\,,
$$
$$
d_j^+=-i\scal{(2V'\frac{d}{dx}+V'')\psi_j,\psi_j}=\int_I (V'\psi_j^2)'=0
$$
(as could be expected for a self-adjoint $H_2$), and
$$
d_j^-=\scal{V''\psi_j,\psi_j}\ge \sqrt{c}=\min_I \sqrt{V''}\,.
$$
Then, by Corollary~\ref{cor:nonsaest1} we have
$$
\min_k |\mu_j-\lam_k|\le \frac{\|V'\|_1^2}{\displaystyle \min_I \sqrt{V''}}\,.
$$
\end{exmp}

\end{document}